\documentclass{birkjour}

\usepackage{color}
 \newtheorem{thm}{Theorem}[section]
 \newtheorem{cor}[thm]{Corollary}
 \newtheorem{lem}[thm]{Lemma}
 \newtheorem{prop}[thm]{Proposition}
 \theoremstyle{definition}
 \newtheorem{defn}[thm]{Definition}
 \theoremstyle{remark}

 \newtheorem*{prob}{Problem}
 \numberwithin{equation}{section}

\renewcommand{\d}{\mathrm{d}}

\begin{document}

%
%
%
%
%
%
%
%
%

\title[Metric dimension of imprimitive graphs]
 {On the metric dimension of imprimitive\\ distance-regular graphs}

\author[R.\ F.\ Bailey]{Robert F.\ Bailey}

\address{%
Division of Science (Mathematics)\\
Grenfell Campus, Memorial University of Newfoundland\\
University Drive\\
Corner Brook, NL A2H~6P9\\
Canada}

\email{rbailey@grenfell.mun.ca}



\date{August 2015}

\begin{abstract}
A {\em resolving set} for a graph $\Gamma$ is a collection of vertices $S$, chosen so that for each vertex $v$, the list of distances from $v$ to the members of $S$ uniquely specifies $v$. The {\em \mbox{metric} dimension} of $\Gamma$ is the smallest size of a resolving set for $\Gamma$.  Much attention has been paid to the metric dimension of distance-regular graphs.  
Work of Babai from the early 1980s yields general bounds on the metric dimension of primitive distance-regular graphs in terms of their parameters.
We show how the metric dimension of an imprimitive distance-regular graph can be related to that of its halved and folded graphs, but also consider infinite families (including Taylor graphs and the incidence graphs of certain symmetric designs) where more precise results are possible.
\end{abstract}

\maketitle

\section{Introduction} \label{section:intro}
Let $\Gamma=(V,E)$ be a finite, undirected graph without loops or multiple edges.  For $u,v\in V$, the {\em distance} from $u$ to $v$ is the least number of edges in a path from $u$ to $v$, and is denoted $\d_\Gamma(u,v)$ (or simply $\d(u,v)$ if $\Gamma$ is clear from the context).

A {\em resolving set} for a graph $\Gamma=(V,E)$ is a set of vertices $R=\{v_1,\ldots,v_k\}$ such that for each vertex $w \in V$, the list of distances $(\d(w,v_1),\ldots,\d(w,v_k))$ uniquely determines $w$.  Equivalently, $R$ is a resolving set for $\Gamma$ if, for any pair of vertices $u,w \in V$, there exists $v_i \in R$ such that $\d(u,v_i)\neq \d(w,v_i)$; we say that $v_i$ {\em resolves} $u$ and $w$.  The {\em metric dimension} of $\Gamma$ is the smallest size of a resolving set for $\Gamma$.  This concept was introduced to the graph theory literature in the 1970s by Harary and Melter~\cite{Harary76} and, independently, Slater~\cite{Slater75}; however, in the context of arbitrary metric spaces, the concept dates back at least as far as the 1950s (see Blumenthal~\cite{Blumenthal53}, for instance).  For further details, the reader is referred to the survey paper \cite{bsmd}.  

When studying metric dimension, distance-regular graphs are a natural class of graphs to consider.  A graph $\Gamma$ with diameter~$d$ is {\em distance-regular} if, for all $i$ with $0\leq i\leq d$ and any vertices $u,v$ with $\d(u,v)=i$, the number of neighbours of $v$ at distances $i-1$, $i$ and $i+1$ from $u$ depend only on the distance $i$, and not on the choices of $u$ and $v$.  
These numbers are denoted by $c_i$, $a_i$ and $b_i$ respectively, and are known as the {\em parameters} of $\Gamma$.  It is easy to see that $c_0$, $b_d$ are undefined, $a_0=0$, $c_1=1$ and $c_i+a_i+b_i=k$ (where $k$ is the valency of $\Gamma$).  We put the parameters into an array, called the {\em intersection array} of $\Gamma$,
\[ \left\{ \begin{array}{cccccc}
\ast & 1   & c_2 & \cdots & c_{d-1} & c_d \\ 
 0   & a_1 & a_2 & \cdots & a_{d-1} & a_d \\
 k   & b_1 & b_2 & \cdots & b_{d-1} & \ast 
\end{array} \right\}. \]
In the case where $\Gamma$ has diameter~$2$, we have a {\em strongly regular} graph, and the intersection array may be determined from the number of vertices~$n$, \linebreak valency~$k$, and the parameters $a=a_1$ and $c=c_2$; in this case, we say $(n,k,a,c)$ are the parameters of the strongly regular graph.  Another important special case of distance-regular graphs are the {\em distance-transitive} graphs, i.e.\ those graphs $\Gamma$ with the property that for any vertices $u,v,u',v'$ such that $\d(u,v)=\d(u',v')$, there exists an automorphism $g$ such that $u^g=u'$ and $v^g=v'$.  For more information about distance-regular graphs, see the book of Brouwer, Cohen and Neumaier~\cite{BCN} and the forthcoming survey paper by van~Dam, Koolen and Tanaka~\cite{vanDam}.  In recent years, a number of papers have been written on the subject of the metric dimension of distance-regular graphs (and on the related problem of class dimension of association schemes), by the present author and others: see \cite{small,jk,bsmd,grassmann,Beardon,Caceres07,FengWang,Gravier,attenuated,dualpolar,GuoWangLi,fourfamilies,Heger,SeboTannier04}, for instance.
In this paper, we shall focus on various classes of {\em imprimitive} distance-regular graphs, which are explained below.

\subsection{Primitive and imprimitive graphs}
A distance regular graph $\Gamma$ with diameter $d$ is {\em primitive} if and only if each of its distance-$i$ graphs (for $0<i\leq d$) is connected, and is {\em imprimitive} otherwise.  For $d=2$, i.e.\ strongly regular graphs, the only imprimitive examples are the complete multipartite graphs $s K_t$ (with $s>1$ parts of size $t$).  For valency $k\geq 3$, a result known as {\em Smith's Theorem} (after D.~H.~Smith, who proved it for the distance-transitive case~\cite{Smith71}) states that there are two ways for a distance-regular graph to be imprimitive: either the graph is bipartite, or is {\em antipodal}.  The latter case arises when the distance-$d$ graph consists of a disjoint union of cliques, so that the relation of being at distance $0$ or $d$ in $\Gamma$ is an equivalence relation on the vertex set.  The vertices of these cliques are referred to as {\em antipodal classes}; if the antipodal classes have size~$t$, then we say $\Gamma$ is {\em $t$-antipodal}.  It is possible for a graph to be both bipartite and antipodal, with the hypercubes providing straightforward examples.

If $\Gamma$ is a bipartite distance-regular graph, the distance-$2$ graph has two connected components; these components are called the {\em halved graphs} of $\Gamma$.  If $\Gamma$ is $t$-antipodal, the {\em folded graph}, denoted $\overline{\Gamma}$, of $\Gamma$ is defined as having the antipodal classes of $\Gamma$ as vertices, with two classes being adjacent in $\overline{\Gamma}$ if and only if they contain adjacent vertices in $\Gamma$.  The folded graph $\overline{\Gamma}$ is also known as an {\em antipodal quotient} of $\Gamma$; conversely, $\Gamma$ is an {\em antipodal $t$-cover} of $\overline{\Gamma}$.  We note that $\Gamma$ and $\overline{\Gamma}$ have equal valency; a result of Gardiner~\cite[Corollary 4.4]{Gardiner74} shows that $t$ is at most this valency.  The operations of halving and folding may be used to reduce imprimitive graphs to primitive ones: see~\cite[{\S}4.2A]{BCN} for details.  In particular, an imprimitive distance-regular graph with valency $k\geq 3$ may be reduced to a primitive one by halving at most once and folding at most once.  Also, the intersection arrays of the halved and/or folded graphs of $\Gamma$ may be obtained from that of $\Gamma$.

In a 2006 paper of Alfuraidan and Hall~\cite[Theorem 2.9]{AH06}, a refinement of Smith's Theorem is obtained which will be especially useful to us.  We summarize their result below.

\begin{thm}[Alfuraidian and Hall] \label{thm:AH}
Let $\Gamma$ be a connected distance-regular graph with $n$ vertices, diameter~$d$ and valency~$k$.  Then one of the following occurs:
\begin{enumerate}
\item $\Gamma$ is primitive, with $d\geq 2$ and $k\geq 3$;
\item $k=2$, and $\Gamma$ is a cycle $C_n$;
\item $d\leq 1$, and $\Gamma$ is a complete graph $K_n$;
\item $d=2$, and $\Gamma$ is a complete multipartite graph $s K_t$ (where $n=st)$;
\item $d=3$, $\Gamma$ is both bipartite and $2$-antipodal, $n=2v$, and $\Gamma$ is $K_{v,v}-I$\linebreak (a complete bipartite graph with a perfect matching deleted);
\item $d=3$, $\Gamma$ is bipartite but not antipodal, $n=2v$, and $\Gamma$ is the incidence graph of a symmetric design with $v$ points and block size $k<v-1$;
\item $d=3$, $\Gamma$ is antipodal but not bipartite, and $\Gamma$ is an antipodal cover of a complete graph $K_{k+1}$;
\item $d=4$, $\Gamma$ is both bipartite and antipodal, $\Gamma$ is an antipodal cover of a complete bipartite graph $K_{k,k}$, and the halved graph is complete multipartite;
\item $d=6$, $\Gamma$ is both bipartite and antipodal, the halved graphs are antipodal of diameter~$3$, and the folded graph is bipartite of diameter~$3$;
\item $d\geq 4$, $\Gamma$ is antipodal but not bipartite, and the folded graph is primitive with diameter $\lfloor d/2 \rfloor$ and valency $k\geq 3$;
\item $d\geq 4$, $\Gamma$ is bipartite but not antipodal, and the halved graphs are primitive with diameter $\lfloor d/2 \rfloor$ and valency at least~$3$;
\item $d$ is odd and $d=2e+1\geq 5$, $\Gamma$ is bipartite and $2$-antipodal, the folded graph is primitive with diameter $e\geq 2$ and valency $k\geq 3$, and the halved graphs are primitive with diameter $e\geq 2$ and valency at least~$3$;
\item $d$ is even and $d=2e\geq 8$, $\Gamma$ is bipartite and antipodal, and the graphs obtained by successive halving and folding are primitive with diameter $\lfloor e/2 \rfloor \geq 2$ and valency at least~$3$.
\end{enumerate}
\end{thm}
We shall refer to these classes as AH1--AH13.  The numbering is not exactly the same as that given in~\cite{AH06}. We have separated the case of $K_{v,v}-I$ from the other imprimitive graphs of diameter~$3$ (cf.~\cite[Corollary 2.10]{AH06}); as we shall see, it is something of a special case.  Apart from $C_3\cong K_3$ and $C_4\cong K_{2,2}$, no graph appears in more than one class.

\subsection{Metric dimension and asymptotics}
One of main aims of this paper is to consider the asymptotic behaviour of the metric dimension of imprimitive distance-regular graphs.  We use Knuth's convention for asymptotic notation (from~\cite{Knuth76}; see also~\cite[{\S}9.2]{concrete}): for $n$ sufficiently large, we say that $f(n)=O(g(n))$ if there is a constant $C$ such that $f(n)\leq C\cdot g(n)$, that $f(n)=\Omega(g(n))$ if there is a constant $C'$ with $f(n)\geq C'\cdot g(n)$, and $f(n)=\Theta(g(n))$ if both of these happen.

For any graph with $n$ vertices and diameter $d$, it is straightforward to see that the metric dimension $\mu$ must satisfy the inequality $n \leq \mu+d^{\mu}$ (see \cite[Proposition 3.6]{bsmd}); when considering families of fixed diameter (such as strongly regular graphs), this gives a lower bound on $\mu$ of $\Omega(\log n)$.  

Like Alfuraidan and Hall~\cite{AH06}, we regard class AH1 as the ``generic'' class of distance-regular graphs.
For graphs in this class, namely primitive distance-regular graphs of diameter $d\geq 2$ and valency $k\geq 3$, the pioneering work of Babai~\cite{Babai80,Babai81} in the early 1980s (in a different context: see \cite[\S{3.4--3.5}]{bsmd} for details) yields the following.

\begin{thm}[Babai]
\label{thm:babai}
Suppose that $\Gamma$ is a primitive distance-regular graph with $n$ vertices, valency $k\geq 3$ and diameter $d\geq 2$.  Then:
  \begin{itemize}
  \item[(i)] $\mu(\Gamma) < 4\sqrt{n}\log n$;
  \item[(ii)] if $d=2$ (i.e.~$\Gamma$ is strongly regular), we have (a) $\mu(\Gamma) < 2\sqrt{n}\log n$ and (b) $\displaystyle \mu(\Gamma)< \frac{2n^2}{k(n-k)}\log n < \frac{4n}{k} \log n$ (where $k\leq n/2$);
  \item[(iii)] if $M(\Gamma)$ is the maximum size of a set of vertices at a given distance from any vertex of $\Gamma$, we have 
  \[ \mu(\Gamma) < 2d \frac{n}{n-M(\Gamma)} \log n. \]
  \end{itemize}
\end{thm}

Some remarks about Theorem~\ref{thm:babai} are in order.  While the bound (b) in part (ii) appears weaker than (a), if $n$ and $k$ have a linear relationship (for example, in the case of Paley graphs where $n=2k+1$) this may be combined with the lower bound above to obtain $\mu(\Gamma)=\Theta(\log n)$.  On the other hand, the strongly regular Johnson and Kneser graphs $J(m,2)$ and $K(m,2)$ have metric dimension $\Theta(\sqrt{n})$ (see~\cite[Corollary 3.33]{bsmd}), so the $\sqrt{n}$ factor cannot be eliminated in general.
In part (iii) the value of $M$ may be calculated from the intersection array; for families where $d$ is fixed we may also obtain an improvement on the $O(\sqrt{n}\log n)$ upper bound.

For classes AH2--AH4, the metric dimension is easy to determine exactly, as we summarise below.

\begin{prop} \label{prop:easy}
For the graphs in families AH2--AH4, we have $\mu(C_n)=2$, $\mu(K_n)=n-1$ and $\mu(s K_t)=s(t-1)$.
\end{prop}

\proof For cycles and complete graphs, this is trivial; for complete multipartite graphs, it is straightforward (see~\cite[Proposition~1]{small}). \endproof

Asymptotically, for both complete graphs and complete multipartite graphs on $n$ vertices, we have $\mu(\Gamma)=\Theta(n)$.  However, we suggest that classes AH2--AH4 are pathological cases which should be disregarded as atypical.  For the remaining classes AH5--AH13, much more work is required.


\section{General results for imprimitive graphs}

Given that imprimitive distance-regular graphs may be reduced to primitive graphs by the operations of halving and folding, it is desirable to obtain relationships between resolving sets and metric dimension for imprimitive graphs and their halved and/or folded graphs.  In this section, we obtain such relationships.

\subsection{Halving and folding}

First, we consider the halving operation when $\Gamma$ is bipartite.  Our first result does not assume that $\Gamma$ is distance-regular.

\begin{thm} \label{thm:halving}
Let $\Gamma=(V,E)$ be a connected bipartite graph with bipartition $V=V^+\cup V^-$, and let $\Gamma^+=(V^+,E^+)$ and $\Gamma^-=(V^-,E^-)$ be its halved graphs.  Then $\mu(\Gamma)\leq\mu(\Gamma^+)+\mu(\Gamma^-)$.
\end{thm}

\proof Let $R^+\subseteq V^+$ and $R^-\subseteq V^-$ be resolving sets for $\Gamma^+$ and $\Gamma^-$, respectively.  We will show that $R=R^+\cup R^-$ is a resolving set for $\Gamma$.

Let $x,y\in V$.  If one of these is in $V^+$ and the other in $V^-$, then for any vertex $w\in V$, we must have that one of $\d_\Gamma(x,w)$ and $\d_\Gamma(y,w)$ is odd and the other is even, and thus $\d_\Gamma(x,w) \neq \d_\Gamma(y,w)$, i.e.\ $w$ resolves $x$ and $y$.  In particular, we are free to choose $w\in R$.

So we must consider the case where $x$ and $y$ are both in the same bipartite half. If $x,y\in V^+$, then there exists $w\in R^+$ such that $\d_{\Gamma^+}(x,w)\neq \d_{\Gamma^+}(y,w)$, and thus
\[ \d_{\Gamma}(x,w) = 2\cdot\d_{\Gamma^+}(x,w)\neq 2\cdot\d_{\Gamma^+}(y,w) = \d_{\Gamma}(y,w), \]
i.e.\ $w$ resolves $x$ and $y$ (in $\Gamma$).  The case where $x,y\in V^-$ is similar.

Hence $R$ is a resolving set for $\Gamma$, and thus $\mu(\Gamma)\leq\mu(\Gamma^+)+\mu(\Gamma^-)$. \endproof

A simple application of this theorem is when $\Gamma=K_{n,n}$, so $\Gamma^+ \cong \Gamma^- \cong K_n$; clearly $\mu(K_{n,n})=2n-2 = 2\mu(K_n)$ and the bound holds with equality.  

Even in the case of distance-regular graphs, the halved graphs need not be isomorphic (although they must have the same parameters), so we cannot assume that $\mu(\Gamma^+)$ and $\mu(\Gamma^-)$ should be equal.  Indeed, a counterexample is provided by the incidence graph of the unique (up to duality) generalized quadrangle $GQ(3,3)$: its halved graphs (i.e.\ the point graphs of the $GQ(3,3)$ and its dual) have metric dimension~$7$ and~$8$ respectively (see~\cite[Table~11]{small}); the incidence graph itself has metric dimension~$10$ (see~\cite[Table~4]{small}).  This also shows that the upper bound in Theorem~\ref{thm:halving} need not be met with equality.

Obtaining bounds on the metric dimension of an antipodal distance-regular graph $\Gamma$ in terms of its folded graph $\overline{\Gamma}$ is less straightforward.  Suppose $\Gamma=(V,E)$ is $t$-antipodal, with $V$ partitioned into $s$ antipodal classes $W_1,\ldots,W_s$.  We define a {\em $t$-antipodal partition} to be a partition of $V$ into transversals of $W_1,\ldots,W_s$.  (Note that some authors use the term ``antipodal partition'' to refer to the partition of $V$ into antipodal classes, so care is required here.)  We remark that a $t$-antipodal graph has many $t$-antipodal partitions, although in certain cases there natural partition arising from how the graph is constructed.

We note that in an antipodal distance-regular graph $\Gamma$ with diameter $d$, for $u,v\in V$ with $\d_{\Gamma}(u,v)=d$ we have $\d_{\Gamma}(u,x)+\d_{\Gamma}(x,v)=d$ for any $x\in V$, and the diameter of $\overline{\Gamma}$ is $\bar{d} = \lfloor d/2 \rfloor$.

\begin{thm} \label{thm:folding}
Let $\Gamma=(V,E)$ be an $(r+1)$-antipodal distance-regular graph with diameter $d$, and let $V=V^0 \cup V^1 \cup \cdots \cup V^r$ be an $(r+1)$-antipodal partition of $V$.  Suppose $\overline{R}$ is a resolving set for the folded graph $\overline{\Gamma}=(\overline{V},\overline{E})$ whose diameter is $\bar{d}$.  For $v\in\overline{V}$, let $\{v^0,v^1,\ldots,v^r\}$ be its inverse image in $V$ with $v^i\in V^i$ for each $i\in\{0,1,\ldots,r\}$.  Let $R=\{w^1,\ldots,w^r\, : \, w\in \overline{R} \}$.
  \begin{itemize}
  \item[(i)] If $u^i,v^j \in V$ with $u,v\in\overline{V}$, $i,j\in\{0,1,\ldots,r\}$ and $u\neq v$, then there exists $w^\ell \in R$ with $w\in \overline{R}$ and $\ell \in \{1,\ldots,r\}$ such that $\d_\Gamma(u^i,w^\ell) \neq \d_\Gamma(v^j,w^\ell)$ (i.e.\ $w^\ell$ resolves $u^i$ and $v^j$).
	
	\item[(ii)] If $d=2\bar{d}+1$ is odd, or for every $u\in\overline{V}$ there exists $w\in\overline{R}$ such that $\d_{\overline{\Gamma}}(u,w)<\bar{d}$, then $R$ is a resolving set for $\Gamma$.  In particular, $\mu(\Gamma) \leq r\cdot \mu(\overline{\Gamma})$.
	
	\item[(iii)] If $d=2\bar{d}$ is even and there exists $u\in\overline{V}$ such that $\d_{\overline{\Gamma}}(u,w)=\bar{d}$ for all $w\in R$, then $R^\ast = R\cup\{u^1,\ldots,u^r\}$ is a resolving set for $\Gamma$.  In particular, $\mu(\Gamma) \leq r(\mu(\overline{\Gamma})+1)$.
	\end{itemize}
\end{thm}

\proof We start by proving part (i) of the theorem.  Suppose that $u^i,v^j \in V$ with $u,v\in\overline{V}$, $i,j\in\{0,1,\ldots,r\}$ and $u\neq v$.  Since $\overline{R}$ is a resolving set for $\overline{\Gamma}$, there exists $w\in\overline{R}$ such that $\d_{\overline{\Gamma}}(u,w) \neq \d_{\overline{\Gamma}}(v,w)$.  Consequently, there exist indices $i',j'$ such that 
\[ \d_{\Gamma}(u^i,w^{i'}) = \d_{\overline{\Gamma}}(u,w) \neq \d_{\overline{\Gamma}}(v,w) = d_{\Gamma}(v^j,w^{j'}). \]
We wish to find some $\ell\in\{1,\ldots,r\}$ so that $w^\ell \in R$ and $w^\ell$ resolves $u^i$ and $v^j$; we consider the following cases.
  \begin{enumerate}
	\item If $i'=j'\geq 1$, let $\ell=i'=j'\in\{1,\ldots,r\}$, so $w^\ell \in R$ and $\d_{\Gamma}(u^i,w^\ell) \neq \d_{\Gamma}(v^j,w^\ell)$.
	
	\item If $i'=j'=0$, we have 
	\[ \d_{\Gamma}(u^i,w^1) = d - \d_{\Gamma}(u^i,w^0) = d - \d_{\overline{\Gamma}}(u,w) \]
	and similarly 
	\[ \d_{\Gamma}(v^j,w^1) = d - \d_{\Gamma}(v^j,w^0) = d - \d_{\overline{\Gamma}}(v,w) \]
	by the properties of antipodal vertices in $\Gamma$.  In particular, we have $\d_{\Gamma}(u^i,w^1) \neq \d_{\Gamma}(v^j,w^1)$, so $w^1\in R$ and resolves $u^i$ and $v^j$.
	
	\item Otherwise, we have $i'\neq j'$. At most one of $i',j'$ can be $0$; without loss of generality we suppose that $i'\neq 0$, and so $w^{i'}\in R$.  For a contradiction, suppose that $w^{i'}$ does not resolve $u^i$ and $v^j$, i.e.\ $\d_{\Gamma}(u^i,w^{i'}) = \d_{\Gamma}(v^j,w^{i'})$.  Since $2\bar{d}\leq d$, we then have
\[ \bar{d} \geq \d_{\overline{\Gamma}}(u,w) = \d_{\Gamma}(u^i,w^{i'}) = \d_{\Gamma}(v^j,w^{i'}) \]
\[ = d-\d_{\Gamma}(v^j,w^{j'}) = d-\d_{\overline{\Gamma}}(v,w) \geq d-\bar{d} \geq \bar{d}. \]
	In particular, this implies that $\bar{d} = \d_{\overline{\Gamma}}(u,w) = \d_{\overline{\Gamma}}(v,w)$, a contradiction.  Thus $w^{i'}$ resolves $u^i$ and $v^j$.
	\end{enumerate}
	
To prove parts (ii) and (iii), we must construct a resolving set for $\Gamma$.  By part (i), we know that for vertices $u^i,v^j\in V$, if $u\neq v$ then there is a vertex in $R$ which resolves them.  It remains to consider the case where $u=v$ and $i\neq j$, i.e.\ $u^i$ and $v^j$ are antipodal vertices in $\Gamma$.

First, we suppose that there exists $w\in \overline{R}$ such that $\d_{\overline{\Gamma}}(u,w)=d'<\bar{d}$.
Then if $i\neq j$, for any $\ell\in\{1,\ldots,r\}$, we must have that 
\[ \{ \d_{\Gamma}(u^i,w^\ell),\d_{\Gamma}(u^j,w^\ell) \} = \{d',d-d'\}. \]
Since $d'<\bar{d}$ and $d\geq2\bar{d}$, we have $d'\neq d-d'$, and thus $w^\ell$ resolves $u^i$ and $u^j$.  
Combined with part (i), this shows that $R$ is a resolving set for $\Gamma$ of size $r\cdot |\overline{R}|$, and thus (ii) holds.

Otherwise, we must have that $\d_{\overline{\Gamma}}(u,w)=\bar{d}$ for all $w\in \overline{R}$.  Since $\overline{R}$ is a resolving set for $\overline{\Gamma}$, there can be at most one vertex of $\overline{\Gamma}$ with this property.

If $d=2\bar{d}+1$ is odd, for any $\ell\in\{1,\ldots,r\}$, we must have that 
\[ \{ \d_{\Gamma}(u^i,w^\ell),\d_{\Gamma}(u^j,w^\ell) \} = \{\bar{d},\bar{d}+1 \} \]
and so $w^\ell$ resolves $u^i$ and $u^j$.  Thus $R$ is a resolving set for $\Gamma$ of size $r\cdot |\overline{R}|$, and again (ii) holds.

If $d=2\bar{d}$ is even, we have $\d_{\Gamma}(u^i,w^\ell)=\d_{\Gamma}(u^j,w^\ell)=\bar{d}$ for all $i,j,\ell$, and thus $u^i,u^j$ are not resolved by any $w^\ell$.  However, by taking all but one of $u^0,\ldots,u^r$ along with $R$, we have that $R^\ast = R\cup\{u^1,\ldots,u^r\}$ is a resolving set for $\Gamma$ of size $r\cdot(|\overline{R}|+1)$, and thus (iii) holds.

This completes the proof. \endproof

Straightforward examples are provided by the complete multipartite graphs $s K_t$: these are $t$-antipodal covers of the complete graph $K_s$, and have diameter~$2$.  Since any resolving set $\overline{R}$ for $K_t$ contains $t-1$ vertices, the remaining vertex is adjacent to all of $\overline{R}$, so case (iii) of Theorem~\ref{thm:folding} applies.  This gives $\mu(s K_t)\leq (t-1)((s-1)+1)=s(t-1)$, which we know from Proposition~\ref{prop:easy} to be the exact value, so the upper bound is achieved.  Further such cases are discussed in Section~\ref{subsection:2antipodal} below.

Some examples of where the bounds in Theorem~\ref{thm:folding} are not achieved can be found in the tables in~\cite{small}.  First, the Conway--Smith graph $E$ on $63$ vertices and with diameter~$4$ is a $3$-antipodal cover of the Kneser graph $K(7,2)$; from~\cite[Table~10]{small} we have $\mu(E)=6$ and $\mu(K(7,2))=4$, so \mbox{$\mu(E)<2\cdot\mu(K(7,2))$.}  Second, the Foster graph $F$ on $90$ vertices and with diameter~$8$ is a $3$-antipodal cover of Tutte's $8$-cage $T$; from~\cite[Table~3]{small} we have $\mu(F)=5$ and $\mu(T)=6$, so $\mu(F)<2\cdot \mu(T)$.  In fact, in this latter case the covering graph has smaller metric dimension than its folded graph.

\subsection{Some consequences} \label{subsection:conseq}

When considering the classification of Theorem~\ref{thm:AH}, the most immediate applications of Theorems~\ref{thm:halving} and~\ref{thm:folding} are to graphs when halving and/or folding yields a primitive graph with diameter at least~$2$; these are precisely the graphs in classes AH10--AH13.  Combined with Babai's Theorem~\ref{thm:babai}, we can obtain upper bounds on the metric dimension of any graph in those classes in terms of its parameters.  First, we have the following result for bipartite graphs in classes AH11--13.

\begin{cor} \label{cor:bipartite-babai}
Suppose $\Gamma$ is a bipartite distance-regular graph with $n$ vertices, valency $k>2$ and diameter $d\geq 4$, and whose halved graphs are primitive and have diameter at least~$2$.  Then we have:
\begin{itemize}
\item[(i)] $\mu(\Gamma) < 4\sqrt{2n}\log (n/2)$;
\item[(ii)] $\displaystyle \mu(\Gamma) < d\frac{n}{n-2M(\Gamma)}\log(n/2)$ (where $M(\Gamma)$ denotes the maximum size of a set of vertices at a given distance from any vertex of $\Gamma$).
\end{itemize}
\end{cor}

\proof Using Theorem~\ref{thm:halving}, we know that $\mu(\Gamma)\leq \mu(\Gamma^+)+\mu(\Gamma^-)$.  Since the halved graphs are primitive distance-regular graphs with diameter $\lfloor d/2 \rfloor$, we can apply Theorem~\ref{thm:babai} to them.  We note that the distance classes of $\Gamma^+$ and $\Gamma^-$ are formed from the distance classes of $\Gamma$, so we have $M(\Gamma)=M(\Gamma^+)=M(\Gamma^-)$.  \endproof

In case AH10, we have that $\Gamma$ is antipodal but not bipartite, so we need to apply Theorem~\ref{thm:folding} instead.

\begin{cor} \label{cor:antipodal-babai}
Suppose $\Gamma$ is a $t$-antipodal distance-regular graph with $n$ vertices, valency $k>2$ and diameter $d\geq 4$, and whose folded graph $\overline{\Gamma}$ is primitive and has diameter at least~$2$.  Then we have:
\begin{itemize}
\item[(i)] $\displaystyle \mu(\Gamma) < (t-1)\left(4\sqrt{n/t}\log(n/t)+1\right)$;
\item[(ii)] $\displaystyle \mu(\Gamma) < (t-1)\left(d\frac{n}{n-tM(\overline{\Gamma})}\log(n/t)+1\right)$ (where $M(\overline{\Gamma})$ denotes the maximum size of a set of vertices at a given distance from any vertex of $\overline{\Gamma}$).
\end{itemize}
\end{cor}

For graphs in classes AH5--AH9, halving and/or folding yields either a complete or complete bipartite graph, and so we obtain upper bounds on the metric dimension of such graphs which are linear in the number of vertices.  Much of the remainder of this paper is devoted to improving upon this.  However, we shall first consider some properties of $2$-antipodal graphs.


\subsection{2-antipodal graphs} \label{subsection:2antipodal}
In the case where $\Gamma$ is $2$-antipodal with diameter $d$ and folded graph $\overline{\Gamma}$, Theorem~\ref{thm:folding} shows that $\mu(\Gamma)\leq \mu(\overline{\Gamma})$ if $d$ is odd, and $\mu(\Gamma)\leq \mu(\overline{\Gamma})+1$ if $d$ is even.  However, we can obtain more detailed results in this case: the following lemma will be especially useful.  For any vertex~$v$ of $\Gamma$, we denote by $\Gamma_i(v)$ the set of vertices of $\Gamma$ that are at distance~$i$ from~$v$.

\begin{lem} \label{lemma:plusminus}
Suppose that $\Gamma$ is a $2$-antipodal distance regular graph of diameter $d$, whose vertex set has a $2$-antipodal partition $V^+\cup V^-$.  Then, without loss of generality, a resolving set for $\Gamma$ can be chosen just from vertices in $V^+$.
\end{lem}

\proof We claim that if $R$ is any resolving set for $\Gamma$ and $v^- \in R$, then \linebreak $(R\setminus \{v^- \}) \cup \{v^+ \}$ is also a resolving set.  To show this, suppose that $x,y$ are resolved by $v^-$, i.e.\ $\d_{\Gamma}(x,v^-)\neq \d_{\Gamma}(y,v^-)$.  Suppose that $\d_{\Gamma}(x,v^+)=i$, i.e.\ $x\in\Gamma_i(v^+)$.  Since $\Gamma$ is distance-regular, there exists a path of length $d-i$ to some vertex in $\Gamma_d(v^+)$; however, as $v^-$ is the unique vertex in $\Gamma_d(v^+)$, it follows that $\d_{\Gamma}(x,v^-)=d-i$, and so $\d_{\Gamma}(x,v^+)+\d_{\Gamma}(x,v^-)=d$.  Therefore, $\d_{\Gamma}(x,v^+)= d-\d_{\Gamma}(x,v^-) \neq d-\d_{\Gamma}(y,v^-) = \d_{\Gamma}(y,v^+)$, and hence $v^+$ also resolves $x,y$.

By repeating the above process as required, an arbitrary resolving set for $\Gamma$ may be transformed into a resolving set consisting only of vertices in $V^+$, and the result follows.  \endproof

If $\Gamma$ is bipartite as well as antipodal, and has odd diameter $d=2e+1\geq 3$, then it is necessarily $2$-antipodal (otherwise, it would contain an odd cycle of length $3d$); such graphs form classes AH5 and AH12.  We have the following theorem.

\begin{thm} \label{thm:bip_2ant_odd}
Let $\Gamma=(V,E)$ be a bipartite, $2$-antipodal distance-regular graph with odd diameter $d=2e+1$, with folded graph $\overline{\Gamma}=(\overline{V},\overline{E})$ of diameter $\bar{d}=e$.  Then $\mu(\Gamma)=\mu(\overline{\Gamma})$.
\end{thm}

\proof  By Theorem~\ref{thm:folding} we have $\mu(\Gamma)\leq \mu(\overline{\Gamma})$.

To show the converse, we consider the bipartition $V=V^+\cup V^-$, which is also a $2$-antipodal partition, so for any $v\in\overline{V}$, its preimages in $V$ are $v^+\in V^+$ and $v^-\in V^-$.
By Lemma~\ref{lemma:plusminus}, there exists a resolving set $R^+$ for $\Gamma$ with $R^+\subseteq V^+$.  Let $\overline{R}=\{ w\in \overline{V} \, : \, w^+ \in R^+ \}$.  For all $u^+,v^+ \in V^+$, there exists $w^+ \in R^+$ with $\d_{\Gamma}(u^+,w^+) \neq \d_{\Gamma}(v^+,w^+)$.  If these distances both lie in the interval $\{0,\ldots,e\}$, we must have
\[ \d_{\overline{\Gamma}}(u,w) = \d_{\Gamma}(u^+,w^+) \neq \d_{\Gamma}(v^+,w^+) = \d_{\overline{\Gamma}}(v,w). \]
Similarly, if these distances both lie in the interval \mbox{$\{e+1,\ldots,d\}$,} we have
\[ \d_{\overline{\Gamma}}(u,w) = d - \d_{\Gamma}(u^+,w^+) \neq d -\d_{\Gamma}(v^+,w^+) = \d_{\overline{\Gamma}}(v,w). \]
Otherwise, we have (without loss of generality) that 
\mbox{$0\leq \d_{\Gamma}(u^+,w^+) \leq e$} and $e+1 \leq \d_{\Gamma}(v^+,w^+) \leq d$, which implies that $\d_{\overline{\Gamma}}(u,w)$ is even and $\d_{\overline{\Gamma}}(v,w)$ is odd.  Therefore, $w$ resolves $u,v\in \overline{V}$, and thus $\overline{R}$ is a resolving set for $\overline{\Gamma}$.  Hence $\mu(\overline{\Gamma}) \leq \mu(\Gamma)$.

This completes the proof. \endproof

Immediately, we have the following corollary about the graphs in class AH5.

\begin{cor} \label{cor:doubledclique}
The metric dimension of the graph $K_{v,v}-I$, i.e.\ a complete bipartite graph with a $1$-factor removed, is $v-1$.
\end{cor}

\proof We have that $K_{v,v}-I$ is a bipartite, $2$-antipodal distance-regular graph with diameter~$d=3$, and its folded graph is the complete graph $K_v$ which has metric dimension $v-1$.  Then we apply Theorem~\ref{thm:bip_2ant_odd}. \endproof

In terms of the asymptotic behaviour of metric dimension, this tells us that for graphs in class AH5, namely $\Gamma=K_{v,v}-I$ with $n=2v$, we have $\mu(\Gamma)=\Theta(n)$, in common with classes AH3 and AH4.  So we may regard this class as another pathological case.  For graphs in class AH12, we have no change in the asymptotics from what we saw in the previous subsection, although we can be more precise.

The following definition gives us an alternative interpretation of Theorem~\ref{thm:bip_2ant_odd}.

\begin{defn} \label{defn:bipdouble}
Let $\Gamma=(V,E)$ be a graph.  Then the {\em bipartite double} (or {\em bipartite cover}) of $\Gamma$ is the bipartite graph $D(\Gamma)$ whose vertex set consists of two disjoint copies of $V$, labelled $V^+$ and $V^-$, and where $u^+ \in V^+$ and $w^- \in V^-$ are adjacent in $D(\Gamma)$ if and only if $u$ and $w$ are adjacent in $\Gamma$.
\end{defn}

For example, the bipartite double of a complete graph $K_{v+1}$ is the graph $K_{v,v}-I$ from Corollary~\ref{cor:doubledclique} above.  More generally, if $\Gamma$ is distance-regular with diameter $d$ and odd girth $2d+1$, then $D(\Gamma)$ is distance-regular with diameter $2d+1$, and is an antipodal $2$-cover of $\Gamma$.  Furthermore, any distance-regular graph of odd diameter which is both bipartite and antipodal must arise this way (see~\cite[{\S}4.2D]{BCN}), and the bipartition $V^+\cup V^-$ is also a $2$-antipodal partition.  Thus Theorem~\ref{thm:bip_2ant_odd} may be rephrased as follows.

\begin{thm} \label{thm:bipdouble}
Suppose that $\Gamma$ is a distance-regular graph of diameter $d$ with odd girth $2d+1$.  Then the metric dimension of its bipartite double $D(\Gamma)$ is equal to the metric dimension of $\Gamma$.
\end{thm}

As an example of its applications, Theorem~\ref{thm:bipdouble} may be applied to the following infinite family.
The {\em Odd graph} $O_k$ has as its vertex set the collection of all 
$(k-1)$-subsets of a $(2k-1)$-set, with two vertices adjacent if and only if the corresponding $(k-1)$-sets are disjoint.  (The Odd graph $O_3$ is the Petersen graph.)  This graph is distance-regular, has diameter $k-1$ and odd girth $2k-1$, so therefore satisfies the conditions of Theorem~\ref{thm:bipdouble}; its bipartite double is known as the {\em doubled Odd graph}.  (See~\cite[{\S}9.1D]{BCN} for further details.)  Consequently, we have another corollary.

\begin{cor} \label{cor:doubledodd}
The Odd graph $O_k$ and doubled Odd graph $D(O_k)$ have equal metric dimension, which is at most $2k-2$.
\end{cor}

\proof It follows from \cite[Theorem 6]{jk} that the Odd graph $O_k$ has metric dimension at most $2k-2$.  Since this graph satisfies the hypotheses of Theorem~\ref{thm:bipdouble}, the result follows.
\endproof

This latter corollary provides a slight improvement on Theorem 3.1 of Guo, Wang and Li \cite{GuoWangLi}, who showed that $\mu(D(O_k)) \leq 2k-1$.  However, it removes the requirement to consider the doubled Odd graph separately from the Odd graph.


\section{Antipodal and diameter 3: Antipodal covers of cliques} \label{section:antipodal3}

For a graph $\Gamma$ in class AH7 of Theorem~\ref{thm:AH}, i.e.\ $\Gamma$ is antipodal, has diameter~$3$ and is not bipartite, the folded graph is a complete graph, and Theorem~\ref{thm:folding} gives a bound on $\mu(\Gamma)$ of $O(n)$.  However, a stronger result is desirable, and is possible in the case where $\Gamma$ is $2$-antipodal: such graphs are called {\em Taylor graphs} and are discussed below.  We shall see that, rather than using the folded graph, we can reduce $\Gamma$ to a primitive strongly regular graph in a different way, which then gives us a suitable relationship to bound the metric dimension.

\subsection{Taylor graphs} \label{subsection:taylor}

A {\em Taylor graph} is a $2$-antipodal distance-regular graph on $2n+2$ vertices, obtained via the following construction, due to Taylor and Levingston~\cite{TaylorLevingston78}.  Suppose that $\Delta=(V,E)$ is a strongly regular graph with parameters $(n,2c,a,c)$.  Construct a new graph $\Gamma$ by taking two copies of the set $V$ labelled as $V^+,V^-$, along with two new vertices $\infty^+,\infty^-$, and defining adjacency as follows: let $\infty^+$ be adjacent to all of $V^+$, $\infty^-$ be adjacent to all of $V^-$, $u^+\sim v^+$ and $u^-\sim v^-$ (in $\Gamma$) if and only if $u\sim v$ (in $\Delta$), and $u^+\sim v^-$ if and only if $u\neq v$ and $u\not\sim v$ (where $\sim$ denotes adjacency).

From the construction, one may verify that $\Gamma$ is indeed distance-regular, $2$-antipodal, and that the folded graph is a complete graph $K_{n+1}$.
The given labelling of the vertices ensures that $v^+$ is the unique antipode of $v^-$, for all $v\in V\cup\{\infty\}$.  
For any vertex $x$ of $\Gamma$, let $\Gamma[x]$ denote the subgraph of $\Gamma$ induced on the set of neighbours of $x$.
The construction ensures that $\Delta$ is isomorphic to both $\Gamma[\infty^+]$ and $\Gamma[\infty^-]$; for any other vertex $x$, $\Gamma[x]$ is also strongly regular with the same parameters, but need not be isomorphic to $\Delta$.  As a simple example, one may use this construction to obtain the icosahedron from a $5$-cycle, which has parameters $(5,2,0,1)$.  For further examples, we refer to the table of strongly regular graphs in Brouwer and Haemers~\cite[{\S}9.9]{BrouwerHaemers2012}.

A {\em two-graph} $\mathcal{D}$ is a pair $(V,\mathcal{B})$, where $V$ is a set and $\mathcal{B}$ is a collection of $3$-subsets of $V$, with the property that any $4$-subset of $V$ contains an even number of members of $\mathcal{B}$.  From any element $x\in V$, one may form a graph with vertex set $V\setminus\{x\}$ by deleting $x$ from all triples which contain it, and taking the resulting pairs as edges; such a graph is a {\em descendant} of $\mathcal{D}$.  The collection of all descendants of $\mathcal{D}$ is referred to as a {\em switching class}, because of the relationship with the operation of Seidel switching; for more information on two-graphs and switching classes, see~\cite{SeidelTaylor78,SpenceHandbook}.

A two-graph is {\em regular} if every $2$-subset of $V$ occurs in a constant number of members of $\mathcal{B}$.  In~\cite{Taylor77}, Taylor proved that the descendants of a regular two-graph on $n+1$ points are necessarily strongly regular graphs with parameters $(n,2c,a,c)$.  Taylor and Levingston~\cite{TaylorLevingston78} subsequently showed the following; see also \cite[{\S}1.5]{BCN} for an account of their work.

\begin{thm}[Taylor and Levingston~\cite{TaylorLevingston78}] \label{thm:taylev} $ $
\begin{itemize}
\item[(i)] An antipodal $2$-cover of $K_{n+1}$ is necessarily a Taylor graph.
\item[(ii)] There exists a one-to-one correspondence between Taylor graphs and regular two-graphs on $n+1$ points.
\item[(iii)] The isomorphism classes of descendants of a regular two-graph $\mathcal{D}$, i.e.\ the members of a switching class of strongly regular graphs with \linebreak parameters $(n,2c,a,c)$, are precisely the isomorphism classes of induced subgraphs $\Gamma[v]$ of the corresponding Taylor graph $\Gamma$.
\end{itemize}
\end{thm}

To confuse matters, the strongly regular graphs arising as the descendants of a regular two-graph associated with the unitary group $\mathrm{PSU}(3,q^2)$, as discovered by Taylor~\cite{Taylor77}, are sometimes referred to as ``Taylor's graph'': see~\cite{Spence92}.  Distance-transitive Taylor graphs were classified in 1992~\cite{Taylor92}.

The main result of this section is to relate the resolving sets for a Taylor graph with those for the descendants of the corresponding regular two-graph.

\begin{thm} \label{thm:taylor}
Let $\mathcal{D}$ be a regular two-graph with corresponding Taylor graph $\Gamma$, and let $\{\Delta_1,\ldots,\Delta_s\}$ be the switching class of descendants of $\mathcal{D}$.  Choose a descendant $\Delta$ with the smallest metric dimension, i.e.\ $\mu(\Delta) \leq \mu(\Delta_i)$ for all descendants $\Delta_i$.  Then we have:
\begin{itemize}
\item[(i)] $\mu(\Gamma) = \mu(\Delta)+1$; and
\item[(ii)] $\mu(\Delta_i) \in \{ \mu(\Delta),\, \mu(\Delta)+1 \}$ for all descendants $\Delta_i$.
\end{itemize}
\end{thm}

\proof  First, we show that $\mu(\Gamma) \leq\mu(\Delta)+1$.  Label the vertices of $\Gamma$ as $V^+ \cup V^- \cup \{\infty^+,\infty^-\}$, as described above, and choose a smallest resolving set $R\subseteq V$ for $\Delta$.

We will show that $R^+\cup\{\infty^+\}$ is a resolving set for $\Gamma$.  Since $R$ is a resolving set for $\Delta$, then for any pair of distinct vertices $u,v\in V$, there exists $x\in R$ such that $\d_\Delta(u,x) \neq \d_\Delta(v,x)$.  Since $\d_\Delta(u,x)=\d_\Gamma(u^+,x^+)$ and $\d_\Delta(v,x)=\d_\Gamma(v^+,x^+)$, it follows that $x^+$ resolves the pair $(u^+,v^+)$.  Likewise, $x^-$ resolves the pair $(u^-,v^-)$; however, since $\Gamma$ is $2$-antipodal, Lemma~\ref{lemma:plusminus} shows that $x^+$ will also resolve the pair $(u^-,v^-)$.  Any pair of vertices of the form $(u^+,v^-)$ will be resolved by $\infty^+$, as $\d_\Gamma(u^+,\infty^+)=1$ for any $u^+\in V^+$, and $\d_\Gamma(v^-,\infty^+)=2$ for any $v^-\in V^-$.  Finally, any pair involving one of $\infty^+$ or $\infty^-$ will be resolved by $\infty^+$, since $\infty^-$ is the unique vertex at distance~$3$ from $\infty^+$.

Now we will establish the reverse inequality, i.e.\ $\mu(\Gamma)\geq \mu(\Delta)+1$.  Choose a resolving set $S$ for $\Gamma$ of size $\mu(\Gamma)$.  Now choose some vertex $x\in S$, and consider the subgraph $\Gamma[x]$ induced on the set $N(x)$ of neighbours of $x$.
Since $\Gamma$ is a Taylor graph, $\Gamma[x]$ must be isomorphic to a descendent $\Delta_i$ of the regular two-graph $\mathcal{D}$, and thus has diameter~$2$.  Furthermore, the vertices in $\{x\}\cup N(x)$ form one part of a $2$-antipodal partition, so by applying Lemma~\ref{lemma:plusminus}, we may assume that the remaining vertices of $S$ are all neighbours of $x$.

Since $S$ is a resolving set for $\Gamma$, then for any $u,v\in N(x)$, there exists a vertex $w\in S$ that resolves the pair $(u,v)$; note that $w\neq x$, as $x$ is clearly adjacent to all of its neighbours.  Furthermore, for any pair of vertices $u,v\in N(x)$, we have that $\d_\Gamma(u,v)=\d_{\Gamma[x]}(u,v)$: since $\Gamma[x]$ is an induced subgraph, $u$ and $v$ are adjacent in $\Gamma$ if and only if they are adjacent in $\Gamma[x]$, while if $u$ and $v$ are not adjacent, they have distance~$2$ in $\Gamma$ (in a path through $x$) and distance~$2$ in $\Gamma[x]$ (since it has diameter~$2$).  As we assumed that $S\setminus\{x\}\subseteq N(x)$, this shows that $S\setminus\{x\}$ is a resolving set of size $\mu(\Gamma)-1$ for $\Gamma[x]$.  
Consequently, we have
\[ \mu(\Delta) \leq \mu(\Delta_i) = \mu(\Gamma[x]) \leq \mu(\Gamma)-1, \]
as required, and this concludes the proof of part (i).

To prove part (ii), we note that a given descendant $\Delta_i$ need not arise in the manner described above, i.e.\ induced on the set of neighbours of a vertex $x$ of a minimum resolving set for $\Gamma$.  However, any resolving set for $\Gamma$ may be used to construct a resolving set of the same size for $\Delta_i$.  Suppose that $\Delta_i\cong \Gamma[w]$ for some vertex $w$.  If $S$ is a minimum resolving set for $\Gamma$ that does not contain $w$, then we can still apply Lemma~\ref{lemma:plusminus} to assume that $S\subseteq N(w)$, and the same argument as above shows that $S$ is also a resolving set for $\Gamma[w]$.  Therefore, $\mu(\Delta_i) \leq \mu(\Gamma)$, and we have
\[ \mu(\Delta) \leq \mu(\Delta_i) \leq \mu(\Gamma) = \mu(\Delta)+1, \]
and part (ii) follows. 
\endproof

We remark that in the case of vertex-transitive Taylor graphs (such as those obtained from Paley graphs), all descendants are isomorphic, and the result simply states $\mu(\Gamma) = \mu(\Delta)+1$.

The result in part (ii) of Theorem~\ref{thm:taylor} seems a little unsatisfactory: a better result would be that all descendants of a given Taylor graph (i.e.\ all strongly \mbox{regular} graphs in the same switching class) have the same metric dimension, although the author was unable to show this.  There is computational evidence to support such a claim.  It is known that strongly \mbox{regular} graphs with the same parameters need not have the same metric \mbox{dimension:} the Paley graph on 29 vertices has metric dimension~$6$, while the other strongly regular graphs with parameters $(29,14,6,7)$, which fall into five switching classes, all have metric dimension~$5$ (see~\cite[Table 2]{small}).  Furthermore, the $3854$ strongly regular graphs with parameters $(35,16,6,8)$, which fall into exactly $227$ switching classes~\cite{McKaySpence01}, all have metric dimension~$6$ (see \linebreak \cite[Table~13]{small}).  (As an application of Theorem~\ref{thm:taylor}, we know that all $227$ Taylor graphs on~$72$ vertices have metric dimension~$7$.)

Given what we know about the metric dimension of primitive strongly regular graphs from Theorem~\ref{thm:babai}, we can combine this with Theorem~\ref{thm:taylor} to obtain bounds on the metric dimension of Taylor graphs.

\begin{cor} \label{cor:taylorbabai}
Suppose that $\Gamma$ is a Taylor graph with $N=2n+2$ vertices.  Then (a) $\mu(\Gamma) < 2\sqrt{n}\log n +1$ and (b) $\displaystyle \mu(\Gamma) < 4\frac{n}{k}\log n +1$ (where $k$ is the valency of a descendant of $\Gamma$).  In particular, $\mu(\Gamma)=O(\sqrt{N}\log N)$.
\end{cor}

Of course, if the descendants of a Taylor graph $\Gamma$ are strongly regular graphs with logarithmic metric dimension, then this carries over to $\Gamma$.  For example, if $\Gamma$ has a Paley graph as a descendant, then this has metric dimension at most $2\log n$ (as shown by Fijav\v{z} and Mohar~\cite{FijavzMohar04}); in this case it follows that $\mu(\Gamma)=\Theta(\log N)$.

\begin{prob}
What can be said about the metric dimension of $t$-antipodal distance-regular graphs of diameter~$3$, where $t>2$?
\end{prob}


\section{Bipartite and diameter 3: Incidence graphs of symmetric designs} \label{section:incidence}

In this section, we will consider graphs in class AH6, namely bipartite distance-regular graphs with diameter~$3$.  Now, if $\Gamma$ is such a graph, its halved graphs of $\Gamma$ will be complete graphs, and so Theorem~\ref{thm:folding} gives an upper bound of $n-2$ on its metric dimension; as in the previous section, it would be desirable to improve on this.  

A {\em symmetric design} (or {\em square $2$-design}) with parameters $(v,k,\lambda)$ is a pair $(X,\mathcal{B})$, where $X$ is a set of $v$ \emph{points}, and $\mathcal{B}$ is a family of $k$-subsets of $X$, called \emph{blocks}, such that any pair of distinct points are contained in exactly $\lambda$ blocks, and that any pair of distinct blocks intersect in exactly $\lambda$ points.  It follows that $|\mathcal{B}|=v$.  A symmetric design with $\lambda=1$ is a {\em projective plane}, while a symmetric design with $\lambda=2$ is known as a {\em biplane}~\cite{Cameron73}.

The {\em incidence graph} of a symmetric design is the bipartite graph with vertex set $X\cup\mathcal{B}$, with the point $x\in X$ adjacent to the block $B\in\mathcal{B}$ if and only if $x\in B$.  It is straightforward to show that the incidence graph of a symmetric design is a bipartite distance-regular graph with diameter~$3$.  The converse is also true (see~\cite[{\S}1.6]{BCN}): any bipartite distance-regular graph of diameter~$3$ gives rise to a symmetric design.  

The {\em dual} of a symmetric design is the design obtained from the incidence graph by reversing the roles of points and blocks; $(X,\mathcal{B})$ and its dual both have the same parameters.  The {\em complement} of a symmetric design $(X,\mathcal{B})$ has the same point set $X$, and block set $\overline{\mathcal{B}} = \{ X\setminus B \, : \, B\in \mathcal{B} \}$.  The incidence graph of $(X,\overline{\mathcal{B}})$ is obtained from that of $(X,\mathcal{B})$ by interchanging edges and non-edges across the bipartition.  If $(X,\mathcal{B})$ has parameters $(v,k,\lambda)$, then $(X,\overline{\mathcal{B}})$ has parameters $(v,v-k,v-2k+\lambda)$.

Suppose $\Gamma$ is the incidence graph of $(X,\mathcal{B})$.  For any distinct points $x,y\in X$ and any distinct blocks $A,B \in \mathcal{B}$, we have $\d_{\Gamma}(x,y) = \d_{\Gamma}(A,B) = 2$, while $\d_{\Gamma}(x,B) = 1$ if $x\in B$ and $\d_{\Gamma}(x,B) = 3$ if $x\not\in B$.  It follows that the incidence graph of a symmetric design and that of its complement have the same metric dimension. Clearly, the incidence graph of a symmetric design and that of its dual also have the same metric dimension, as these graphs are isomorphic.

We observe that if a resolving set $R$ for $\Gamma$ is contained entirely within $X$ or entirely within $\mathcal{B}$, we have $|R|\geq v-1$.  In the case of the unique symmetric design with $k=v-1$, where the blocks are all the $(v-1)$-subsets, the graph obtained is $K_{v,v}-I$, which has metric dimension $v-1$ by Corollary~\ref{cor:doubledclique}.  From now on, we shall assume that $k<v-1$, where it is natural to ask if smaller resolving sets exist which must therefore contain both types of vertex.  One approach for constructing resolving sets is as follows.

Suppose $\Gamma$ is the incidence graph of a symmetric design $(X,\mathcal{B})$.  A {\em split resolving set} for $\Gamma$ is a set $R=R_X\cup R_{\mathcal{B}}$, where $R_X\subseteq X$ and $R_\mathcal{B} \subseteq \mathcal{B}$, chosen so that any two points $x,y$ are resolved by a vertex in $R_\mathcal{B}$, and any two blocks $A,B$ are resolved by a vertex in $R_X$.  We call $R_X$ and $R_\mathcal{B}$ {\em semi-resolving sets} for the blocks and points of the design.  The smallest size of a split resolving set will be denoted by $\mu^\ast(\Gamma)$.  We note that a split resolving set is itself a resolving set, as any vertex will resolve a pair $x,B$, given that the parities of the distances to $x$ and to $B$ will be different; therefore, we only need consider resolving point/block pairs.  Clearly, we have $\mu(\Gamma)\leq \mu^\ast(\Gamma)$.

A straightforward observation is that, for any two blocks $A,B$, the point $x$ resolves the blocks $A,B$ if and only if $x$ lies in exactly one of the two blocks (i.e.\ $x\in A$ and $x\not\in B$, or vice-versa), and a block $B$ resolves the points $x,y$ if and only if exactly one of $x,y$ lies in $B$.

\subsection{Projective planes}
In the case of projective planes, the blocks of the design are usually referred to as {\em lines}, and are denoted by $\mathcal{L}$.  It is known that for a projective plane to exist, we have $v=q^2+q+1$ and $k=q+1$ for some integer~$q$, called the {\em order} of the projective plane.  We let $\Gamma_\Pi$ denote the incidence graph of a projective plane $\Pi$.

A {\em blocking set} for a projective plane $\Pi=(P,\mathcal{L})$ of order $q$ is a subset of points $S\subseteq P$ chosen so that every line $L\in\mathcal{L}$ contains at least one point in $S$; moreover, $S$ is a {\em double blocking set} if every line $L$ contains at least two points in $S$.  Ball and Blokhuis~\cite{BallBlokhuis96} showed that, for $q>3$, a double blocking set has size at least $2(q+\sqrt{q}+1)$, with equality occurring in the plane $\mathrm{PG}(2,q)$ when $q$ is a square.  Also, one can easily construct a double blocking set of size $3q$ by taking the points of three non-concurrent lines.  Double blocking sets and semi-resolving sets are related by the following straightforward proposition.

\begin{prop} \label{prop:blocking}
A double blocking set with a single point removed forms a semi-resolving set for the lines of a projective plane.
\end{prop}

\proof Let $S$ be a double blocking set for $\Pi=(P,\mathcal{L})$.  Any pair of distinct lines $L_1,L_2$ intersects in a unique point $x$.  Since $S$ is a double blocking set, there exists $y\in L_1\setminus \{x\}$ such that $y\in S$ and $y\not\in L_2$.  Hence $y$ resolves the lines $L_1,L_2$.  By the same argument, there also exists $z\in L_2\setminus\{x\}$ such that $z\in S$ and $z\not\in L_1$.  This redundancy allows us to delete any point $x$ from $S$ so that $S\setminus\{x\}$ is still a semi-resolving set; however, deleting two points from $S$ may prevent us from resolving some pairs of lines.
\endproof

Using Proposition~\ref{prop:blocking}, we may obtain a split resolving set for $\Gamma_\Pi$ of size $(\tau_2(\Pi)-1)+(\tau_2(\Pi^\perp)-1)$ by taking a semi-resolving set of this form for the points along with the dual of such a set for the lines (where $\tau_2(\Pi)$ denotes the smallest size of a double blocking set in $\Pi$, and $\Pi^\perp$ denotes the dual plane). If $\Pi$~is self-dual then this simplifies as $2(\tau_2(\Pi)-1)$.  At the problem session of the 2011 British Combinatorial Conference, the author asked whether this was best possible.  In 2012, the question was answered by H\'eger and Tak\'ats~\cite{Heger} for the Desarguesian plane $\mathrm{PG}(2,q)$.

\begin{thm}[%
{H\'eger and Tak\'ats %
    \cite[Theorem 4]{Heger}%
    }%
    ] \label{thm:ht1}
A semi-resolving set for $\mathrm{PG}(2,q)$ has size at least $\min\{ 2q+q/4-3,\, \tau_2(\mathrm{PG}(2,q))-2 \}$; for a square prime power $q\geq 121$, this is at least $q+\sqrt{q}$.
\end{thm}

Of course, a minimum resolving set for $\Gamma_\Pi$ need not be a split resolving set.  W.~J.~Martin (personal comunication) was able to construct a non-split resolving set for $\Gamma_\Pi$ of size $4q-4$ (see~\cite[Figure~1]{Heger}), and conjectured that this was best possible (except for small orders).  This conjecture was also proved in the 2012 paper of H\'eger and Tak\'ats~\cite{Heger}.

\begin{thm}[%
{H\'eger and Tak\'ats %
    \cite[Theorem 2]{Heger}%
    }%
    ] \label{thm:ht2}
For a projective plane $\Pi$ of order $q\geq 23$, the metric dimension of its incidence graph $\Gamma_\Pi$ is $\mu(\Gamma_\Pi)= 4q-4$.
\end{thm}

H\'eger and Tak\'ats also gave a complete description of all resolving sets of this size: see~\cite[{\S}3]{Heger}.  Asymptotically, their result gives the following.

\begin{cor} \label{cor:PGasymp}
For a projective plane $\Pi$ whose incidence graph $\Gamma_\Pi$ has $N$ vertices, we have $\mu(\Gamma_\Pi)=\Theta(\sqrt{N})$.
\end{cor}

\proof For a projective plane $\Pi$ of order $q$, we know that the number of vertices of $\Gamma_\Pi$ is $N=2(q^2+q+1)$, and by Theorem~\ref{thm:ht2} (for $q\geq 23$), the metric dimension is $\mu(\Gamma_\Pi)=4q-4=\Theta(\sqrt{N})$.  \endproof

One might ask if this result holds for symmetric designs with $\lambda>1$.  In the next subsection, we consider this possibility.

\subsection{Symmetric designs with a null polarity}
A {\em polarity} of a symmetric design $(X,\mathcal{B})$ is a bijection $\sigma\, : \, X\to \mathcal{B}$ which preserves the point/block incidence relation.  It is straightforward to see that $(X,\mathcal{B})$ admits a polarity if and only if there is an ordering of the points and blocks so that the incidence matrix of the design is symmetric.  A point is called {\em absolute} if it is incident with its image under $\sigma$.  A polarity $\sigma$ is said to be {\em null} if no points are absolute.\footnote{Sometimes, the term ``null polarity'' is used when {\em all} points are absolute; however, this is equivalent to the complement of the design having no absolute points.}  In this situation, the incidence matrix has zero diagonal, and so is the adjacency matrix of a graph $\Delta$; this graph is strongly regular with parameters $(v,k,\lambda,\lambda)$ (see~\cite[Theorem~2.1]{Rudvalis71}).  We observe that a symmetric design $(X,\mathcal{B})$ may admit more than one null polarity, and the corresponding strongly regular graphs need not be isomorphic.  The reader is referred to the book of Ionin and Shrikhande~\cite[{\S}7.4]{Ionin06} for more details, and for several constructions of infinite families of such designs, in particular families arising from Hadamard matrices.

Conversely, if one has a strongly regular graph $\Delta=(V,E)$ with parameters $(v,k,\lambda,\lambda)$, the bipartite double of that graph (recall Definition~\ref{defn:bipdouble}) is the incidence graph of a symmetric design with parameters $(v,k,\lambda)$, which \mbox{admits} a null polarity in an obvious way: the points and blocks may be \mbox{labelled} by $V^+$ and $V^-$ respectively, and the map $\sigma\, : \, v^+ \mapsto v^-$ is a null polarity.  We note that non-isomorphic graphs may give rise to the same symmetric design: for instance, the $4\times 4$ lattice $H(2,4)$ and the Shrikhande graph are non-isomorphic strongly regular graphs with parameters $(16,6,2,2)$, yet their bipartite doubles are isomorphic (and give rise to a $(16,6,2)$-biplane).

Given this relationship with bipartite doubles, one may ask if there is a result similar to Theorem~\ref{thm:bipdouble} which can be applied here to find the metric dimension of $\Gamma$, and we have the following theorem.

\begin{thm} \label{thm:nullpolarity}
Let $\Gamma$ be the incidence graph of a non-trivial $(v,k,\lambda)$ symmetric design with a null polarity, and let $\Delta$ be a corresponding strongly regular graph with parameters $(v,k,\lambda,\lambda)$.  Then $\mu(\Gamma) \leq 2\mu(\Delta)$.
\end{thm}

\proof Since $\Gamma$ is the bipartite double of $\Delta$, we will label the points and blocks of the design by $V^+$ and $V^-$ respectively.  Suppose $R\subseteq V$ is a resolving set for $\Delta$; we will show that $R^+\cup R^-$ is a resolving set for $\Gamma$.  Now, for any distinct $u,w\in V$, we have $\d_{\Gamma}(u^+,w^+) = \d_{\Gamma}(u^-,w^-) = 2$, while $\d_{\Gamma}(u^+,w^-) = 1$ if $u\sim w$ in $\Delta$ and $\d_{\Gamma}(u^+,w^-) = 3$ if not.
Clearly, any vertex resolves $u^+,w^-$ (as the distances will have different parities), so it suffices to consider resolving pairs of vertices of the form $u^+,w^+$ and $u^-,w^-$.

If $u\in R$, then clearly $u^+$ resolves the pair $u^+,w^+$ and $u^-$ resolves the pair $u^-,w^-$ (and likewise if $w\in R$), so we assume that $u,w\not\in R$.  Since $R$ is a resolving set for $\Delta$, there exists $x\in R$ where $\d_\Delta(u,x)\neq \d_\Delta(w,x)$; without loss of generality, this implies that $u\sim x$ and $w\not\sim x$, so therefore $\d_\Gamma(u^+,x^-)=1$ and $\d_\Gamma(w^+,x^-)=3$, and thus $x^-$ resolves the pair $u^+,w^+$.  Similarly, $x^+$ resolves the pair $u^-,w^-$.  Hence any pair of vertices of $\Gamma$ is resolved by a vertex in $R^+\cup R^-$, and we are done.
\endproof

Immediately, we have the following corollary, which is reminiscent of Corollary~\ref{cor:taylorbabai} for Taylor graphs.

\begin{cor} \label{cor:symmbabai}
Let $\Gamma$ be the incidence graph of a non-trivial $(v,k,\lambda)$ symmetric design where $k\leq v/2$ and which has a null polarity.  Then we have (a) $\mu(\Gamma) < 4\sqrt{v}\log v$ and (b) $\displaystyle \mu(\Gamma)< \frac{4v^2}{k(v-k)}\log v < \frac{8v}{k} \log v$.
\end{cor}

\proof Let $\Delta$ be the strongly regular graph associated with the design (as in Theorem~\ref{thm:nullpolarity}).  Theorem~\ref{thm:babai} gives us upper bounds on $\mu(\Delta)$, and the result follows.
\endproof

Since $\Gamma$ has $N=2v$ vertices, Corollary~\ref{cor:symmbabai} gives $\mu(\Gamma) = O(\sqrt{N}\log N)$.  
In cases where $v$ and $k$ have a linear relationship, Babai's stronger bound yields an upper bound of $O(\log v)$ on $\mu(\Delta)$.  There are a number of infinite families of such designs, in particular arising from Hadamard matrices (see~\cite[{\S}7.4]{Ionin06}, and also~\cite{HaemersXiang}); the incidence graphs of such designs therefore have metric dimension $\Theta(\log N)$.

We remark, however, that Theorem~\ref{thm:nullpolarity} can never be applied to the incidence graphs of projective planes, since it is known that a projective plane cannot admit a null polarity (see \cite[Proposition 4.10.1]{BrouwerHaemers2012}).

\begin{prob}
What happens if we remove the hypothesis that the design has a null polarity?  Is it still true that the incidence graph $\Gamma$ of a non-trivial symmetric design with $\lambda>1$ has metric dimension $O(\sqrt{N} \log N)$?  Do the bounds of Corollary~\ref{cor:symmbabai} still hold?
\end{prob}

While the method used in the proof of Theorem~\ref{thm:nullpolarity} depends on the existence of the null polarity (to exploit the relationship with strongly regular graphs), there are instances of parameter sets $(v,k,\lambda)$ with multiple isomorphism classes of designs, of which only some have a null polarity, yet all of their incidence graphs have the same metric dimension.  For example, there are five non-isomorphic symmetric  designs with parameters $(15,8,4)$, giving rise to four non-isomorphic \mbox{incidence} graphs; only one of these designs (the complementary design of the projective geometry $\mathrm{PG}(3,2)$) has a null polarity, yet all four graphs have metric dimension~$8$ (see~\cite[Table 2]{small}).  Also, there are three non-isomorphic $(16,6,2)$-biplanes, with three non-isomorphic incidence graphs; again, only one has a null polarity, yet all three have metric dimension~$8$ (see~\cite[Table 2]{small}).


\section{Conclusion}

In this paper, we have begun the systematic study of the metric dimension of imprimitive distance-regular graphs.  Using Theorems~\ref{thm:halving} and~\ref{thm:folding}, combined with the results of Babai (Theorem~\ref{thm:babai}) for primitive graphs, we saw that bounds on the metric dimension of graphs in classes AH9--12 (of Theorem~\ref{thm:AH}) follow immediately.  The difficult cases are bipartite or antipodal graphs of diameter~$3$ (whose halved or folded graphs are complete) or graphs of diameter~$4$ which are both bipartite and antipodal (where the folded graph is complete bipartite).  In Sections~\ref{section:antipodal3} and~\ref{section:incidence} we obtained some results on these classes, where the diameter is~$3$; obtaining general bounds for the metric dimension of graphs in these classes remains open.  It seems plausible that the bounds in Corollary~\ref{cor:symmbabai} should hold for symmetric designs in general, if the requirement on the existence of a null polarity is removed.  The other cases (namely antipodal covers of complete and complete bipartite graphs) also have connections with design theory (see~\cite{AldredGodsil,GodsilHensel}) and this relationship may prove to be useful when studying these cases.

\newpage

\section*{Acknowledgment}
The author would like to thank the anonymous referees whose suggestions led to significant improvements to this paper.
This research was begun while the author was a postdoctoral fellow at the University of Regina and Ryerson University.  The author would like to thank the Faculty of Science at each university for financial support.

\end{document}